\def\ver{Oct. 5, 2002, v.6}
\documentstyle{amsppt}
\magnification=1200
\hsize=6.5truein
\vsize=8.9truein
\hoffset=-10pt
\voffset=-30pt
\topmatter
\title Nonvanishing of External Products\\
for Higher Chow Groups
\endtitle
\author Andreas Rosenschon and Morihiko Saito
\endauthor
\abstract Consider an external product of a higher cycle and a usual
cycle which is algebraically equivalent to zero.
Assume there exists an algebraically closed subfield
$ k $ such that the higher cycle and its ambient variety are defined
over
$ k $, but the image of the usual cycle by the Abel-Jacobi map is not.
Then we prove that the external product is nonzero if the image of
the higher cycle by the cycle map to the reduced Deligne cohomology
does not vanish.
We also give examples of indecomposable higher cycles on even
dimensional hypersurfaces of degree at least four in a projective
space which satisfy the last condition.
\endabstract
\endtopmatter
\tolerance=1000
\baselineskip=12pt

\def\ssbul{\raise.2ex\hbox{${\scriptscriptstyle\bullet}$}}

\def\bC{{\Bbb C}}
\def\bD{{\Bbb D}}
\def\bP{{\Bbb P}}
\def\bQ{{\Bbb Q}}
\def\bR{{\Bbb R}}
\def\bZ{{\Bbb Z}}
\def\cA{{\Cal A}}
\def\cB{{\Cal B}}
\def\cD{{\Cal D}}
\def\cL{{\Cal L}}
\def\cO{{\Cal O}}
\def\cV{{\Cal V}}
\def\tet{\widetilde{\eta}}
\def\tom{\widetilde{\omega}}
\def\tf{\widetilde{f}}
\def\tH{\widetilde{H}}
\def\tX{\widetilde{X}}
\def\tY{\widetilde{Y}}
\def\tZ{\widetilde{Z}}
\def\tu{\widetilde{u}}
\def\tDe{\widetilde{\Delta}}
\def\tet{\widetilde{\eta}}
\def\of{\overline{f}}
\def\og{\overline{g}}
\def\oh{\overline{h}}

\def\oQ{\overline{\Bbb Q}}

\def\an{\text{\rm an}}
\def\alg{\text{\rm alg}}
\def\BM{\text{\rm BM}}
\def\ind{\text{\rm ind}}
\def\hom{\text{\rm hom}}
\def\CH{\hbox{{\rm CH}}}
\def\Hdg{\hbox{{\rm Hdg}}}
\def\Ker{\hbox{{\rm Ker}}}
\def\Pic{\hbox{{\rm Pic}}}
\def\Res{\hbox{{\rm Res}}}
\def\div{\text{{\rm div}}\,}
\def\Sing{\text{{\rm Sing}}\,}
\def\Spec{\text{{\rm Spec}}\,}
\def\Ext{\hbox{{\rm Ext}}}
\def\Gr{\text{{\rm Gr}}}
\def\Hom{\hbox{{\rm Hom}}}
\def\MHS{\text{{\rm MHS}}}

\def\SameAuthor{\vrule height3pt depth-2.5pt width1cm}

\document
\centerline{\bf Introduction}

\bigskip\noindent
Let
$ X, Y $ be smooth complex projective varieties.
Take a higher cycle
$ \zeta \in \CH^{p}(X,m) $ [8] and a usual cycle
$ \eta \in \CH^{q}(Y) $, and consider the product
$ \zeta \times \eta \in \CH^{p+q}(X\times Y,m) $.
We assume
$ m > 0 $ and
$ \eta $ is homologically equivalent to zero.
Then the image of
$ \zeta \times \eta $ in
$ \bQ $-Deligne cohomology vanishes, and we are interested in
the question:
When is
$ \zeta \times \eta $ nonzero in
$ \CH^{p+q}(X\times Y,m)_{\bQ} $?
Note that if the varieties and the cycles are defined over number
fields the product is expected to vanish as a consequence of
conjectures of Beilinson [2], [4] and Bloch [7].
We obtained a partial answer when
$ X, Y $ and one of
$ \zeta, \eta $ are defined over an algebraically closed
subfield
$ k $ of
$ \bC $,
but the remaining cycle is not defined over
$ k $ in an appropriate sense, see [32].
In general, this is a rather difficult problem.
A similar problem was studied by C. Schoen in the case of usual
cycles, see [40], [41].
In this paper we consider the case where
$ X, \zeta $ are defined over an algebraically closed subfield
$ k $, but
$ Y, \eta $ are not.
Then we get a variation of Hodge structure on the base
space of a model of
$ Y $, and the problem becomes much more difficult.
We assume
$ \eta $ is algebraically equivalent to zero, and the image of
$ \eta $ by the Abel-Jacobi map is not defined over
$ k $ (more precisely, the image of
$ \eta $ does not come from a
$ k $-valued point of an abelian variety over
$ k $ whose base change by
$ k \to \bC $ is an abelian subvariety of
$ J_{\alg}^{q}(Y) $;
note that there is a largest abelian subvariety defined over
$ k $, which is called the
$ K/k $-trace, see [27]).

\medskip\noindent
{\bf 0.1.~Theorem.}
{\it Let
$ \zeta $,
$ \eta $ be as above {\rm (}i.e.
$ \eta $ is algebraically equivalent to zero,
$ X $,
$ \zeta $ are defined over
$ k $, but
$ Y $ and the image of
$ \eta $ by the Abel-Jacobi map are not{\rm ).}
Assume further that the image of
$ \zeta $ by the cycle map to the reduced
$ \bQ $-Deligne cohomology does not vanish.
Then
$ \zeta \times \eta \ne 0 $ in
$ \CH^{p+q}(X\times Y,m)_{\bQ} $.
}

\medskip
Here the reduced Deligne cohomology is the usual Deligne cohomology if
$ m > 1 $, and is the quotient of it by the group of Hodge cycles
tensored with
$ \bC^{*} $ if
$ m = 1 $, see [31].
If the algebraic part of the Jacobian has a model as an abelian scheme
over a proper
$ k $-variety, we can prove the assertion of Theorem (0.1)
for the cycle map to the usual Deligne cohomology.
If we do not assume that
$ \eta $ is algebraically equivalent to zero, we would have to
assume that
$ m \ge 2q $, and the above condition on
$ \eta $ should be replaced by the cohomological condition (2.5.2)
on the cycle class, see (2.6).
The proof of Theorem (0.1) uses the ``spreading out'' of algebraic
cycles (see [6], [22], [43], [44]) together with the cycle map to
Deligne cohomology ([3], [9], [19], [20], [21], [24], [36]).
The arguments are similar to those in the theory of refined
cycle maps associated to arithmetic mixed sheaves or Hodge
structures, see [32], [37], [38] (and also [1]).
We also use the theory of mixed Hodge modules to get a bound of
weights of the cohomology of a variation of Hodge structure,
see the proof of (2.2).
Actually the arguments show the nonvanishing of
$ \zeta \times \eta $ in
$ \Gr_{F_{L}}^{2}\CH^{p+q}(X\times Y,m)_{\bQ} $,
where
$ F_{L} $ is induced from the Leray filtration by using the
refined cycle map, see (2.7).

As for the conditions in Theorem (0.1), it is easy to construct
an example of an elliptic curve which is not defined over a given
algebraically closed subfield
$ k $ of
$ \bC $ (hence the
$ \bC/k $-trace is trivial) by using the
$ j $-invariant.
Concerning the hypothesis on the higher cycle,
we show that there are many examples satisfying the
hypothesis in the case
$ m = 1 $.

\medskip\noindent
{\bf 0.2.~Theorem.}
{\it For any positive integers
$ n, d $ such that
$ d \ge 4 $, there are smooth hypersurfaces
$ X $ of degree
$ d $ in
$ \bP^{2n+1} $ together with a higher cycle
$ \zeta \in \CH^{n+1}(X,1) $ whose image by the cycle map to
the reduced
$ \bQ $-Deligne cohomology does not vanish.
In particular,
$ \zeta $ is a nontrivial indecomposable cycle.
}

\medskip
See Theorems (3.3) and (4.4).
Note that an indecomposable higher cycle on an odd dimensional
hypersurface cannot be detected by the reduced cycle map,
because the cohomology of a hypersurface is essentially trivial
except for the middle degree.
Theorem (0.2) generalizes results of Collino [13] and of
del Angel and M\"uller-Stach [14] in the case of quartic surfaces.
For another example satisfying the last hypothesis of (0.1),
see [36].
Examples of indecomposable higher cycles on hypersurfaces of degree
$ 2n $ in
$ \bP^{n+1} \,(n\ge 2) $ are constructed by Voisin [45];
however, for
$ n > 2 $, these cycles cannot be detected by the reduced cycle map.
For other examples of indecomposable cycles, see [12], [31], etc.
It does not seem that examples of indecomposable higher cycles
on hypersurfaces of arbitrarily high degree have been known in
the literature.
The proof of (3.3) uses the theory of degeneration of Hodge
structures and period integrals ([11], [39], [42], etc.)
It actually shows that the transcendental part of the image of the
indecomposable cycle by the cycle map does not vanish
(see [14], [36] for other such examples).
The proof of (4.4) was inspired by the Thom-Sebastian theorem for
vanishing cycles.

The paper is organized as follows.
In Sect.~1, we review some basic facts from the theory of the
cycle map of higher Chow groups to Deligne cohomology,
and also the
$ K/k $-trace of an abelian variety.
Then we prove Theorem (0.1) in Sect.~2.
Examples of indecomposable higher cycles on surfaces are
constructed in Sect.~3,
and the higher dimensional case is treated in Sect.~4.

\bigskip\bigskip
\centerline{{\bf 1.~Preliminaries}}

\bigskip\noindent
{\bf 1.1.~Cycle map to Deligne cohomology.}
Let
$ X $ be a smooth complex algebraic variety, and
$ m $ a nonnegative integer.
Let
$ \CH^{p}(X,m)_{\bQ} $ be the higher Chow group with
rational coefficients [8].
We have a cycle map ([9], [19], [36])
$$
cl : \CH^{p}(X,m)_{\bQ} \to H_{\cD}^{2p-m}(X,\bQ(p)),
\leqno(1.1.1)
$$
where the target denotes the absolute Hodge cohomology [3].
The latter coincides with
$ \bQ $-Deligne cohomology ([20], [21], [24]) if
$ X $ is smooth proper.
We have a short exact sequence
$$
\aligned
0 \to \Ext_{\MHS}^{1}(\bQ,H^{2p-m-1}(X,\bQ)
&(p)) \to H_{\cD}^{2p-m}(X,\bQ(p))
\\
&\to\Hom_{\MHS}(\bQ,H^{2p-m}(X,\bQ)(p)) \to 0.
\endaligned
\leqno(1.1.2)
$$

Assume
$ X $ is smooth proper.
Then for
$ m > 0 $, (1.1.2) gives an isomorphism
$$
H_{\cD}^{2p-m}(X,\bQ(p)) = \Ext_{\MHS}^{1}(\bQ,H^{2p-m-1}
(X,\bQ)(p)).
\leqno(1.1.3)
$$
If
$ m = 1 $, we can naturally identify
$ \Hdg^{p-1}(X)_{\bQ}\otimes_{\bZ}\bC^{*} $ with a subspace of (1.1.3),
because
$ \Ext_{\MHS}^{1}(\bQ,\bQ(1)) = \bC^{*}\otimes_{\bZ}\bQ $.
Here\
$ \Hdg^{p-1}(X) $ denotes the group of Hodge cycles of codimension
$ p - 1 $.
We define the reduced Deligne cohomology to be
$$
\Ext_{\MHS}^{1}(\bQ,H^{2p-2}(X,\bQ(p))/
\Hdg^{p-1}(X)_{\bQ}\otimes_{\bZ}\bC^{*}.
\leqno(1.1.4)
$$

If
$ m = 0 $,
the composition of the cycle map (1.1.1) with the last morphism
of (1.1.2) is the usual cycle map, and (1.1.1) induces the
Abel-Jacobi map to the intermediate Jacobian [23] (tensored with
$ \bQ) $:
$$
\CH_{\hom}^{p}(X)_{\bQ} \to J^{p}(X)_{\bQ} =
\Ext_{\MHS}^{1}(\bQ,H^{2p-1}(X,\bQ)(p)).
\leqno(1.1.5)
$$
Here
$ \CH_{\hom}^{p}(X) $ is the subgroup consisting of cycles
homologically equivalent to zero, and the last isomorphism
follows from [10].

\medskip\noindent
{\bf 1.2.~Algebraic part of intermediate Jacobian.}
For a smooth complex projective variety
$ X $, let
$ \CH_{\alg}^{p}(X) $ denote the subgroup consisting of
cycles algebraically equivalent to zero, and let
$ J_{\alg}^{p}(X) $ denote its image in
$ J^{p}(X) $.
Then
$ J_{\alg}^{p}(X) $ is an abelian subvariety of
$ J^{p}(X) $,
and is the image of
$$
\Pic^{0}(\tZ) \to J^{p}(X)
\leqno(1.2.1)
$$
for the normalization
$ \tZ $ of some closed subvariety
$ Z $ of pure codimension
$ p - 1 $ in
$ X $,
see [38], 3.10.
So there exists a finitely generated subfield
$ K $ of
$ \bC $ such that
$ \tZ $ and hence
$ \Pic^{0}(\tZ) $ are defined over
$ K $.
We can verify that the image of (1.2.1) is also defined over
$ K $ by Chow's theorem (see [27]), replacing
$ K $ with a finite extension if necessary.

\medskip\noindent
{\bf 1.3. K/k-trace of abelian variety.}
Let
$ \cA $ be an abelian variety defined over a field
$ K $ of characteristic
$ 0 $, and
$ k $ be an algebraically closed subfield of
$ K $.
Then there exists a largest abelian subvariety
$ \cB $ which is defined over
$ k $ by [27] (because
$ k $ is algebraically closed).
It is called the
$ K/k $-trace of
$ \cA $.
This is well-defined because
$$
\Hom_{k}(\cB_{1},\cB_{2})
= \Hom_{K}(\cB_{1}\otimes_{k}K,\cB_{2}\otimes_{k}K)
$$
for abelian varieties
$ \cB_{1}, \cB_{2} $ over
$ k $ by Chow's theorem (loc.~cit).
We say that a
$ K $-valued point of
$ \cA $ is defined over
$ k $,
if it comes from a
$ k $-valued point of the
$ K/k $-trace of
$ \cA $.

Let
$ L $ be a field containing
$ K $,
and set
$ \cA_{L} = \cA\otimes_{K}L $.
Then the
$ L/k $-trace of
$ \cA_{L} $ does not necessarily coincide with the
$ K/k $-trace of
$ \cA $.
But they coincide if we replace
$ K $ with a finite extension and
$ \cA $ with the base change (using again Chow's theorem).

\bigskip\bigskip
\centerline{{\bf 2.~Proof of Main Theorem}}

\bigskip\noindent
{\bf 2.1. Abelian scheme and variation of Hodge structure.}
Let
$ \cA $ be an abelian scheme over a smooth complex algebraic
variety
$ S $, and
$ H $ be the corresponding polarizable variation of
$ \bZ $-Hodge structure of weight
$ -1 $ and level
$ 1 $,
see [16].
Then we have canonical injective morphisms (see [35])
$$
\cA(S) \to \Ext^{1}(\bZ,H) \to \cA^{\an}(S^{\an}),
\leqno(2.1.1)
$$
where the first term is the group of algebraic sections of
$ \cA \to S $,
the extension group in the middle term is taken in the
category of admissible variations of mixed Hodge structures, and
the last term is the group of analytic sections.
(The first injection is an isomorphism at least if
$ \dim S = 1 $.)
We have furthermore a short exact sequence
$$
0 \to \Ext_{\MHS}^{1}(\bZ,H^{0}(S,H)) \to \Ext^{1}(\bZ,H)
\to \Hom_{\MHS}(\bZ,H^{1}(S,H)) \to 0.
\leqno(2.1.2)
$$
(See [34] and also [46] for the case
$ \dim S = 1 $.)
Note that the canonical injective morphism
$ a_{X}^{*}H^{0}(S,H) \to H $ is a morphism of variations of
Hodge structures [16], where
$ a_{X} : X \to \Spec \bC $ denotes the structure morphism.
In particular,
$ H^{0}(S,H) $ is a polarizable Hodge structure of level
$ 1 $, and the pull-back of the corresponding abelian variety
gives the maximal abelian subscheme of
$ \cA $ coming from an abelian variety on
$ \Spec \bC $.

For a dense open subvariety
$ U $ of
$ S $, (2.1.2) implies the injectivity of
$$
\Hom_{\MHS}(\bZ,H^{1}(S,H)) \to
\Hom_{\MHS}(\bZ,H^{1}(U,H)),
\leqno(2.1.3)
$$
using the snake lemma together with the injectivity of the last
morphism of (2.1.1).

\medskip\noindent
{\bf 2.2.~Proposition.}
{\it Let
$ H_{1} $ be a
$ \bQ $-Hodge structure of weight
$ r $,
and
$ \tH_{2} $ be a polarizable variation of
$ \bQ $-Hodge structure of weight
$ -1 $ and level
$ 1 $ associated with an abelian scheme over an irreducible
smooth complex algebraic variety
$ S $.
Let
$ \xi_{1} \in \Ext_{\MHS}^{1}(\bQ,H_{1}) $,
$ \xi_{2} \in \Hom_{\MHS}(\bQ,H^{1}(S,\tH_{2})) $.
Assume they are nonzero,
$ r \le -2 $ and
$ \Hom_{\MHS}(\bQ(1),H_{1}) = 0 $.
Then for any nonempty open subvariety
$ U $ of
$ S $,
the image of
$ \xi_{1}\otimes \xi_{2} $ by the canonical morphism
$$
\Ext_{\MHS}^{1}(\bQ,H_{1}\otimes H^{1}(S,\tH_{2})) \to
\Ext_{\MHS}^{1}(\bQ,H_{1}\otimes H^{1}(U,\tH_{2}))
$$
does not vanish.
}

\medskip\noindent
{\it Proof.}
Let
$ H_{2} = H^{1}(U,\tH_{2}) $.
Then it has weights
$ \ge 0 $ by construction, see [33].
Furthermore it has weights
$ \le 2 $ and
$ \Gr_{2}^{W}H_{2} $ is of type
$ (1,1) $ (i.e. isomorphic to a direct sum of
$ \bQ(-1) $),
because the underlying Hodge filtration
$ F $ of
$ H_{2} $ satisfies
$ \Gr_{F}^{p} = 0 $ for
$ p > 1 $.
(Indeed, the Hodge filtration comes from the Hodge filtration on
the logarithmic de Rham complex, see loc.~cit).
Consider the long exact sequence
associated with the cohomological functor
$ \Ext^{i}(\bQ,*) $ applied to
$$
0 \to H_{1}\otimes \Gr_{0}^{W}H_{2} \to H_{1}\otimes
H_{2} \to H_{1}\otimes (H_{2}/W_{0}H_{2}) \to 0.
$$
Then the above assertion on the weights and type of
$ H_{2} $ implies the injectivity of
$$
\Ext_{\MHS}^{1}(\bQ,H_{1}\otimes \Gr_{0}^{W}H_{2}) \to
\Ext_{\MHS}^{1}(\bQ,H_{1}\otimes H_{2}).
\leqno(2.2.1)
$$
Since
$ \xi_{2} $ comes from a morphism
$ \xi'_{2} : \bQ \to \Gr_{0}^{W}H_{2} $,
which splits by semisimplicity, and does not vanish by shrinking
$ S $ (see (2.1.3)), the assertion follows.

\medskip\noindent
{\bf 2.3.~Remark.}
The assertion holds for
$ r = -2 $ without assuming the vanishing of
$ \Hom_{\MHS}(\bQ(1),H_{1}) $, if the local monodromies around the
divisor at infinity are semisimple.
Indeed, we can show that
$ H_{2} $ has weights
$ \le 1 $ in this case.

\medskip\noindent
{\bf 2.4.~Lemma.} {\it
Let
$ L $ be a field containing an algebraically closed field
$ k $.
Let
$ \cA $ be an abelian scheme defined over a
$ k $-variety
$ S_{k} $,
and
$ \cA_{L} $ be its pull-back to
$ S_{L} := S_{k}\otimes_{k}L $.
Let
$ \sigma $ be a section of
$ \cA \to S_{k} $, and
$ \sigma_{L} $ be its base change by
$ k \to L $.
If there exists an abelian variety
$ \cB_{L} $ on
$ \Spec L $ together with an injective morphism of abelian
schemes
$ \cB_{L}\times_{L}S_{L}\to \cA_{L} $ such that
$ \sigma_{L} $ comes from a section
$ \sigma'_{L} $ of
$ \cB_{L} $, then
$ \sigma $ satisfies a similar property {\rm (}with
$ L $ replaced by
$ k) $.
}

\medskip\noindent
{\it Proof.}
If such
$ \cB_{L} $ and
$ \sigma'_{L} $ exist, we may assume
$ L $ is finitely generated over
$ k $,
and
$ \cA_{L} $,
$ \cB_{L} $,
$ \sigma'_{L} $ and the morphism are defined over a finitely
generated
$ k $-subalgebra
$ R $ of
$ L $.
Then it is enough to restrict to the fiber over a closed point of
$ \Spec R $.

\medskip\noindent
{\bf 2.5.~Proof of Theorem (0.1).}
By hypothesis, there exist an algebraically closed subfield
$ k $, smooth
$ k $-varieties
$ X_{k} $,
$ \tY_{k}, S_{k} $ together with a proper smooth morphism
$ f : \tY_{k} \to S_{k} $ such that
$ S_{k} $ is irreducible,
$ X = X_{k}\otimes_{k}\bC, $ the base change of the
generic fiber
$ \tY_{K} $ of
$ f $ by an embedding
$ K := k(S_{k}) \to \bC $ is isomorphic to
$ Y $, and
$ \zeta $,
$ \eta $ come from
$ \zeta_{k} \in \CH^{p}(X_{k},m) $,
$ \tet_{k} \in \CH^{q}(\tY_{k}) $.
Furthermore, we have an abelian scheme
$ \cA $ over
$ S_{k} $ such that the base change of the generic fiber of
$ \cA $ by
$ K \to \bC $ is the algebraic part of the intermediate
Jacobian
$ J^{q}(Y) $ and the image of
$ \eta $ by the Abel-Jacobi map is identified with a section
of
$ \cA \to S_{k} $ (replacing
$ S_{k} $ if necessary), because we may assume that
$ \eta $ comes from
$ \Pic^{0}(\tZ) $ in (1.2.1) so that
$ \tet_{k} $ defines a section of
$ \cA $.
We may also assume that the
$ K/k $-trace of the generic fiber of
$ \cA $ coincides with the
$ \bC/k $-trace of the base change by replacing
$ S_{k} $ if necessary, see (1.3).

Let
$ \cA_{\bC} $ be the pull-back of
$ \cA $ to
$ S := S_{k}\otimes_{k}\bC $,
and
$ \tH_{2} $ be the corresponding variation of
$ \bQ $-Hodge structure of weight
$ -1 $ and level
$ 1 $ on
$ S $,
see [16].
This is a direct factor of
$ \tH'_{2} := R^{2q-1}f_{*}\bQ_{\tY}(q) $,
where
$ \tY = \tY_{k}\otimes_{k}\bC $.
We have a natural morphism (see (2.1))
$$
\cA_{\bC}(S)_{\bQ} \to \Hom_{\MHS}(\bQ,H^{1}(S,\tH_{2}))
\subset \Hom_{\MHS}(\bQ,H^{1}(S,\tH'_{2})).
\leqno(2.5.1)
$$
Let
$ \xi_{2} \in \Hom_{\MHS}(\bQ,H^{1}(S,\tH_{2})) $
denote the image of
$ \tet := \tet_{k}\otimes_{k}\bC $.
Then
$$
\xi_{2} \ne 0.
\leqno(2.5.2)
$$
Indeed, if it vanishes, we see that the image of a multiple of
$ \tet $ in
$ \cA_{\bC}(S) $ comes from a section of an abelian variety,
see (2.1).
This contradicts the hypothesis by (2.4), and (2.5.2)
follows.

Let
$ H_{1} = H^{2p-m-1}(X,\bQ)(p) $,
and
$ \xi_{1} \in \Ext_{\MHS}^{1}(\bQ,H_{1}) $ denote the image of
$ \zeta $ by the cycle map (1.1.1).
Then we get
$$
\xi_{1}\otimes \xi_{2} \in
\Ext_{\MHS}^{1}(\bQ,H_{1}\otimes H^{1}(S,\tH_{2})),
$$
which is the image of
$ \zeta \times \tet \in \CH^{p+q}(X\times \tY,m) $ by the
cycle map (1.1.1).
If
$ \zeta \times \eta $ vanishes in
$ \CH^{p+q}(X\times Y,m)_{\bQ} $, there is a dominant morphism
$ \pi $ of a smooth variety
$ S'_{k} $ to
$ S_{k} $ such that the base change of
$ \zeta_{k} \times \tet_{k} $ by
$ \pi $ vanishes, because
$ \zeta \times \eta $ is the base change of
$ \zeta_{k} \times \tet_{k} $ by
$ \Spec \bC \to S_{k} $.
Thus it is enough to show that
$ \xi_{1}\otimes \xi_{2} $ does not vanish after replacing
$ S_{k} $ with any smooth variety
$ S'_{k} $ having a dominant morphism to
$ S_{k} $ (and replacing the cycle by the base change).
But we have a closed subvariety of
$ S'_{k} $ which is finite \'etale over
$ S_{k} $ by shrinking
$ S_{k} $ and
$ S'_{k} $ if necessary.
So we may assume that
$ S'_{k} $ is an open subvariety of
$ S_{k}, $ and the assertion follows from (2.2).
This completes the proof of (0.1).

\medskip\noindent
{\bf 2.6.~Remark.}
If we do not assume that
$ \eta $ is algebraically equivalent to zero but only homologically
equivalent to zero, then the assertion holds for
$ m \ge 2q $ if we assume the condition (2.5.2).
This follows from an estimate of weights of the cohomology of an
algebraic variety in [16].

\medskip\noindent
{\bf 2.7.~Theorem.}
{\it Let
$ F_{L} $ be the filtration of the higher Chow group induced from
the cycle map associated with the theory of arithmetic mixed sheaves
{\rm [37].}
Then with the assumptions of Theorem {\rm (0.1),} the external product
$ \zeta\times\eta $ belongs to
$ F_{L}^{2}\CH^{p+q}(X\times Y,m)_{\bQ} $, and is nonzero in
$ \Gr_{F_{L}}^{2}\CH^{p+q}(X\times Y,m)_{\bQ} $.
}

\medskip\noindent
{\it Proof.}
This follows from (2.5) by using the forgetful
functor from the category of mixed sheaves to that of mixed Hodge
Modules.

\medskip\noindent
{\bf 2.8.~Remark.}
The category of mixed sheaves is a natural generalization of
that of systems of realizations which consist of Betti, de Rham and
$ l $-adic realizations, see [17], [18], [25].
In our situation, we assume that
$ k $ is an algebraically closed subfield of
$ \bC $.
So
$ l $-adic sheaves [5] are not necessary, and we get the category
of mixed Hodge structures whose
$ \bC $-part
$ (H_{\bC};F,W) $ has a
$ k $-structure, i.e.
a bifiltered
$ k $-vector space
$ (H_{k};F,W) $ together with an isomorphism
$ (H_{\bC};F,W) = (H_{k};F,W)\otimes_{k}\bC $ is given.
Thus we get the category of arithmetic mixed Hodge structures
(see also [1]).
Note that the
$ k $-structure has not been used in the above proof.
Forgetting about the
$ k $-structure, we get a variant which is similar to a
formulation of M. Green in the case where the variety is defined
over
$ k $.
(His theory was explained in his talk at Alg. Geom. 2000 Azumino,
Nagano.)
Assuming the conjecture of Beilinson and Bloch on the injectivity
of Abel-Jacobi map for cycles defined over number fields,
we expect to get still the same filtration
$ F_{L} $ on the higher Chow groups after the above modification
in the case
$ k = \oQ $.
It is further expected that this filtration
$ F_{L} $ would give the conjectural filtration of Beilinson [4]
and Bloch [6], see also [26].

\bigskip\bigskip
\centerline{{\bf 3.~Examples of higher cycles on surfaces}}

\bigskip\noindent
{\bf 3.1.~Higher Abel-Jacobi map.}
Let
$ X $ be a smooth proper complex algebraic variety.
A higher cycle
$ \zeta \in \CH^{p}(X,1) $ is represented by
$ \sum_{j} (Z_{j},g_{j}) $ where the
$ Z_{j} $ are irreducible closed subvarieties of codimension
$ p - 1 $ in
$ X $ and the
$ g_{j} $ are rational functions on
$ Z_{j} $ such that
$ \sum_{j} \div g_{j} = 0 $ as a cycle on
$ X $ (without any equivalence relations), see e.g. [31].
We say that
$ \zeta $ is decomposable if the
$ g_{j} $ are constant.
Let
$ \CH_{\ind}^{p}(X,1)_{\bQ} $ be the quotient group of
$ \CH^{p}(X,1)_{\bQ} $ by the subgroup of decomposable
cycles.
An element of
$ \CH_{\ind}^{p}(X,1)_{\bQ} $ is called an indecomposable higher
cycle.
Note that the cycle map induces a well-defined map of
$ \CH_{\ind}^{p}(X,1)_{\bQ} $ to (1.1.4), because the image of
a decomposable cycle is contained in
$ \Hdg^{p-1}(X)_{\bQ}\otimes_{\bZ}\bC^{*} $.

For a higher cycle
$ \zeta = \sum_{j} (Z_{j},g_{j}) $,
let
$ \gamma_{j} $ be the closure of the pull-back of the open
interval
$ (0,+\infty ) $ by
$ g_{j} $.
Then
$ \gamma := \sum_{j} \gamma_{j} $ is a topological cycle of
dimension
$ 2d + 1 $,
where
$ d = \dim X - p $.
It vanishes in
$ H_{2d+1}(X,\bQ)(-d) \,(= H^{2p-1}(X,\bQ)(p)) $,
because it gives the cycle class of
$ \zeta $ in
$ \Hom_{\MHS}(\bQ,H^{2p-1}(X,\bQ)(p)) $ which vanishes
by a weight argument.
So there exists a
$ C^{\infty} $-chain
$ \Gamma $ with
$ \bQ $-coefficients on
$ X $ such that
$ \partial \Gamma = \gamma $.

By Carlson [10], the extension group
$ \Ext_{\MHS}^{1}(\bQ,H^{2p-2}(X,\bQ)(p)) $ is isomorphic to
$$
H^{2p-2}(X,\bC)/\bigl(H^{2p-2}(X,\bQ)(p) +
F^{p}H^{2p-2}(X,\bC)\bigr).
\leqno(3.1.1)
$$
Then the cycle class
$ cl(\zeta) $ in (3.1.1) is represented by a current
$ \Phi_{\zeta} $ defined by
$$
\Phi_{\zeta}(\omega) = (2\pi i)^{-d-1}\biggl(\sum_{j}
\int_{Z_{j}\setminus \gamma_{j}} (\log g_{j})\omega +
2\pi i\int_{\Gamma}\omega\biggr),
\leqno(3.1.2)
$$
where
$ \omega $ is a closed
$ C^{\infty} $-form of type
$ \{(d+1,d+1), \dots, (2d+2,0)\} $.
This formula is due to Beilinson ([3], pp.~61-62) in the case of
$ \bR $-Deligne cohomology, and it is generalized by Levine [28] to
the case of
$ \bZ $-Deligne cohomology
(this construction coincides with the usual
definition of the cycle map, see [24], [36]).

In the case
$ p = \dim X = 2 $ and
$ d = 0 $, consider the quotient of (3.1.1)
$$
H^{2}(X,\bC)/\bigl(H^{2}(X,\bQ)(2) + F^{1}H^{2}(X,\bC)\bigr).
\leqno(3.1.3)
$$
This is called the transcendental part of (3.1.1).
The image of
$ \Hdg^{1}(X)_{\bQ}\otimes_{\bZ}\bC^{*} $ in (3.1.3) vanishes
because the Hodge cycle classes are contained in
$ F^{1}H^{2}(X,\bC) $.
Furthermore, the cycle class
$ cl(\zeta) $ in (3.1.3) is given by the integration
$ \int_{\Gamma}\omega $ for holomorphic
$ 2 $-forms
$ \omega $, i.e. the first term of (3.1.2) vanishes.

\medskip\noindent
{\bf 3.2.~Construction.}
Let
$ (z_{0}, z_{1}, z_{2}, z_{3}) $ be homogeneous coordinates of
$ \bP^{3} $, and
$ d $ an integer such that
$ d \ge 4 $.
For
$ 0 \le i < d $, let
$ f_{i}, g, h \in \bC[z_{1}, z_{2}, z_{3}] $ be homogeneous
polynomials of degree
$ i, 3 $ and
$ d - 3 $ respectively.
Let
$$
f = \sum_{i=0}^{d} z_{0}^{d-i}f_{i}\quad
\text{with}\,\,f_{d} = gh,
$$
and
$ X = f^{ -1}(0) \subset \bP^{3} $,
$ Z = g^{-1}(0) \subset \bP^{2} $.
Consider the parameter space
$ S $ of
$ f_{i}, g, h $ for
$ 0 \le i < d $ such that
$ X $ is smooth and
$ Z $ is a divisor with normal crossings which is either a union of
three lines in
$ \bP^{2} $ or a union of two nonsingular rational curves intersecting
at two points or a rational curve with one ordinary double point.
More precisely, we assume that
$ g $ is one of the following:
$$
z_{1}z_{2}(z_{1}+z_{2}-\lambda z_{3}),\,\,\,
(z_{1}-\lambda^{2} z_{3})(z_{1}z_{3}-z_{2}^{2}),\,\,\,
z_{1}^{2}(z_{1}+\lambda^{2} z_{3})-z_{2}^{2}z_{3},
\leqno(3.2.1)
$$
where
$ \lambda \in \bC^{*} $ is generic.
(In the third case, this gives an example of an indecomposable higher
cycle whose support is irreducible, see [36] for another such example.)
As in the proof of Theorem (3.3) below, we may restrict
$ S $ to a subspace such that
$ f_{i} = 0 $ for
$ 0 < i < d - 1 $, and
$ f_{d-1}, h $ are fixed polynomials.
(Note that a generic member of
$ S $ does not give a generic surface in
$ \bP^{3} $, because the Picard number of
$ X = f^{-1}(0) $ as above is not
$ 1 $.)

We see easily that
$ X $ is smooth near the intersection with
$ z_{0}^{-1}(0) $ if and only if
$$
f_{d-1}^{-1}(0) \cap \Sing f_{d}^{-1}(0) = \emptyset.
\leqno(3.2.2)
$$
In this case,
$ X $ has at most isolated singularities, and
$ X $ is smooth if
$ f_{0} \in \bC $ is general.
Indeed, it is enough to assume that
$ - f_{0} $ is different from a critical value of the function
$ \sum_{0<i\le d} f_{i}/z_{0}^{i} $ on
$ \{z_{0} \ne 0\} $.

Taking rational functions on the irreducible components of
$ Z $ as in (3.1), we get a higher cycle
$ \zeta $ in
$ \CH_{\ind}^{2}(X,1) $ up to a sign, and this gives a family of
higher cycles (modulo decomposable cycles) parametrized by
$ \lambda $ in (3.2.1).
Here we assume that the rational functions have only simple zeros
and poles.
These functions are unique up to constant
multiples and inverses, because there is a unique rational function on
$ \bP^{1} $ up to a constant multiple, which has a simple zero at the
origin and a simple pole at infinity (and a rational curve with one
ordinary double point is obtained by identifying the origin and the
point at infinity of
$ \bP^{1} $).
Note that the family is parametrized a priori by a double cover of the
$ \lambda $-plane which may be non rational, but we need that it is
parametrized by a rational curve to show that the normal function
$ \{\tau_{g}\} $ is constant in the proof of (3.3).

For
$ s \in S $,
let
$ X_{s}, Z_{s}, \zeta_{s} $ denote the corresponding
$ X, Z, \zeta $.
Consider the algebraic family
$ \{X_{s}\} $ over
$ S $.
It gives a holomorphic vector bundle
$ \cV := \{H^{2}(X_{s},\cO_{X_{s}})\}_{s\in S} $ and a local system
$ L := \{H_{2}(X_{s},\bQ)\}_{s\in S} $ over
$ S $.
A multivalued section
$ \eta $ of
$ L $ determines a multivalued holomorphic section
$ \sigma_{\eta} $ of
$ \cV $ by the integrals
$ \int_{\eta}\omega_{i} $, where
$ \{\omega_{i}\} $ is a local basis of the vector bundle
$ \{\Gamma (X_{s},\Omega_{X_{s}}^{2})\}_{s\in S} $.
Let
$ \sigma_{\zeta} $ denote a holomorphic (local) section of
$ \cV $ which is defined by the integrals
$ \int_{\Gamma}\omega_{i} $ as in (3.1).
Although
$ \sigma_{\zeta} $ depends on the choice of
$ \Gamma $,
the ambiguity comes from
$ \sigma_{\eta} $ for some local section
$ \eta $ of
$ L $,
and
$ \sigma_{\zeta} $ gives the cycle map to (3.1.3).
Let
$$
\Sigma = \{s \in S \,| \,\sigma_{\zeta}(s) =
\sigma_{\eta}(s)\,\,\, \text{for some}\,\,\eta \in L \}.
$$

\medskip\noindent
{\bf 3.3.~Theorem.}
{\it With the above notation and assumptions, we have
$ \Sigma \ne S $.
Hence
$ \Sigma $ is locally a countable union of proper analytic closed
subvarieties, and for
$ s \notin \Sigma $,
the cycle class of
$ \zeta $ in {\rm (3.1.3)} and {\rm (1.1.4)} does not vanish;
in particular,
$ \zeta $ is a nontrivial indecomposable cycle.
}

\medskip\noindent
{\it Proof.}
By analytic continuation, it is enough to show that
$ \Sigma \ne S $,
and we may restrict to any open subset of
$ S $ (in the classical topology).
We will derive a contradiction by assuming
$ \Sigma = S $.
We assume that
$ f_{i} = 0 $ for
$ 0 < i < d - 1 $, and choose and fix reduced real polynomials
$ f_{d-1}, h $ such that
$ f_{d-1}^{-1}(0) $ is smooth and intersects
$ h^{-1}(0) $ transversally at smooth points.
Let
$ S' $ be a rational curve with coordinate
$ \lambda $ which parametrizes
$ g $ as in (3.2.1).
Restricting
$ S' $ to an open subvariety, we assume that if
$ g $ belongs to
$ S' $, then
$ f_{d-1}^{-1}(0) $ intersects
$ f_{d}^{-1}(0) $ transversally at smooth points (i.e.
$ g^{-1}(0) $ intersects
$ f_{d-1}^{-1}(0) $ transversely at smooth points, and does not meet
$ f_{d-1}^{-1}(0) \cap h^{-1}(0) $).
Then, replacing
$ S $ with a subvariety, we may assume that
$ S $ is an open subvariety of
$ \bC\times S' $ which is the parameter space of
$ f_{0}, g $ (where
$ f_{i} = 0 $ for
$ 0 < i < d - 1 $, and
$ f_{d-1}, h $ are the fixed polynomials as above).

Let
$ y_{i} = z_{i}/z_{3} $ be affine coordinates of
$ \{z_{3} \ne 0\} \subset \bP^{3} $.
Let
$ \of = f/z_{3}^{d} $,
$ \of_{i} = f_{i}/z_{3}^{i} $,
$ \og = g/z_{3}^{3} $,
$ \oh = h/z_{3}^{d-3} $ so that
$ \of = y_{0}^{d}\of_{0} + y_{0}\of_{d-1} + \of_{d} $ and
$ \of_{d} = \og\,\oh $.
Let
$ \omega_{s} $ be a global
$ 2 $-form on
$ X_{s} $ whose restriction to
$ \{z_{3} \ne 0\} $ is the residue of
$ dy_{0}\wedge dy_{1}\wedge dy_{2}/\of $ along
$ X_{s} $.
This form has a zero of order
$ d - 4 $ along
$ X_{s} \cap \{z_{3} = 0\} $.
Restricted to the open subset which is \'etale over the
$ (y_{1},y_{2}) $-plane,
$ \omega_{s} $ is given by
$ dy_{1}\wedge dy_{2}/(\partial \of/\partial y_{0}) $.
We may assume that
$ f_{d-1} $ does not vanish at
$ O' := \{ z_{1} = z_{2} = 0 \} $, replacing the coordinates
if necessary.
Then
$ \partial \of/\partial y_{0} \ne 0 $ at
$ (y_{0},y_{1},y_{2}) = 0 $.

Let
$ g_{0} \in S' $ corresponding to
$ \lambda = 0 $ in (3.2.1).
Then
$ g_{0} $ is a real polynomial,
$ \Sing Z = \{O'\} $, and
$ Z $ is either three lines meeting at one point,
or two smooth rational curves tangenting at one point,
or a rational curve with one cusp.
Let
$ \Delta_{\varepsilon} $ denote an open disk of radius
$ \varepsilon $.
If
$ \varepsilon, \varepsilon' $ are sufficiently small, and
$ s = (f_{0},g) $ is sufficiently close to
$ (0,g_{0}) $, then
$$
(\partial \of/\partial y_{0})(y_{0},y_{1},y_{2}) \ne 0\,\,\,
\text{for}\,\,(y_{0},y_{1},y_{2}) \in \Delta_{\varepsilon}^{2}
\times \Delta_{\varepsilon'}.
\leqno(3.3.1)
$$
Let
$ \pi $ be the projection of
$ X_{s} \cap \{z_{0} \ne 0\} $ to the
$ (y_{1},y_{2}) $-plane, where
$ X_{s} = f^{-1}(0) $ with
$ f $ as above.
Then, replacing
$ \varepsilon, \varepsilon' $ if necessary, we may assume for
$ s = (f_{0},g) $ sufficiently close to
$ (0,g_{0}) $
$$
\pi : X_{s} \cap \Delta_{\varepsilon}^{2}\times
\Delta_{\varepsilon'} \to \Delta_{\varepsilon}^{2}\,\,\,
\text{is an isomorphism}.
\leqno(3.3.2)
$$

Take
$ s = (f_{0},g) \in \bC\times S' $ sufficiently close to
$ (0,g_{0}) $ such that
$ g $ is real.
Then the higher cycles
$ \zeta_{s} $ and the real
$ 2 $-chains
$ \Gamma_{s} $ on
$ X_{s} $ are defined as above.
More precisely,
$ \Gamma_{s} $ is \'etale by (3.3.2) over an area
$ \Gamma'_{s} $ in
$ \bR^{2} $ which is a connected component of
$ \{\og > 0 \} $ or
$ \{\og < 0 \} $, and is surrounded by the
$ 1 $-chain
$ \gamma_{s} $ contained in
$ Z_{s} \cap \bR^{2} $.
(Here we assume that the rational functions on
$ Z_{s} $ are also defined over
$ \bR $.)
Furthermore, if
$ f_{0} $ is real,
$ \int_{\Gamma}\omega_{s} $ does not vanish by (3.3.1)
because
$ dy_{1}\wedge dy_{2}/(\partial \of/\partial y_{0}) $ is a real
$ 2 $-form and
$ \partial \of/\partial y_{0} $ does not vanish as a function on
$ \Gamma'_{s} $ using (3.3.2).
In particular,
$ \int_{\Gamma}\omega_{s} $ is a nonconstant function of
$ g $ with
$ f_{0} $ fixed, because it vanishes at
$ g_{0} $.
This holds also for
$ f_{0} = 0 $.

We fix
$ g $ as above for a moment, and consider
$ \int_{\Gamma}\omega_{s} $ as a function of
$ f_{0} $ (i.e. we identify
$ f_{0} $ with
$ s) $.
Let
$ \Delta^{*} $ be a sufficiently small punctured disk such that
$ \Delta^{*}\times \{g\} \subset S $.
Then for a multivalued section
$ \eta $ of
$ L $,
$ \sigma_{\zeta} $ and
$ \sigma_{\eta} $ are defined as (multivalued) sections of
$ \cV $ over
$ \Delta^{*}\times \{g\} $.
Since we assume
$ \Sigma = S $, we have
$$
\sigma_{\zeta} = \sigma_{\eta}\quad \text{on}\,\,\Delta^{*}
\times \{g\}\,\,\, \text{for some}\,\,\eta \in L.
\leqno(3.3.3)
$$
In particular,
$ \int_{\Gamma}\omega_{s} = \int_{\eta}\omega_{s} $ on
$ \Delta^{*}\times \{g\} $.
We will derive a contradiction by showing that the limit of
$ \int_{\eta}\omega_{s} $ for
$ s \to (0,g) $ is a constant function of
$ g $, but the corresponding limit of
$ \int_{\Gamma}\omega_{s} $ is not.

By (3.3.3),
$ \sigma_{\eta} $ is a univalent section.
Let
$ \cL $ be the Deligne extension of
$ L|_{\Delta^{*}\times \{g\}} $ over
$ \Delta\times \{g\} $, see [15].
Then the Hodge filtration
$ F $ is extended to
$ \cL $ by [39], and
$ \Gr_{F}^{0}\cL|_{\Delta^{*}\times\{g\}} =
\cV|_{\Delta^{*}\times\{g\}} $.
By Lemma (3.5) below,
$ \sigma_{\eta} $ is extended to a section of
$ \Gr_{F}^{0}\cL $ and its image in
$ \Gr_{F}^{0}\Gr_{V}^{0}\cL(0) $ coincides with the image of some
$ \eta' \in \Gamma (\Delta^{*},L) $
(see (3.4) for the filtration
$ V $.)
By Poincar\'e duality,
$ \eta' $ is identified with a section of the local system
$ \{H^{2}(X_{s},\bQ)\} $.
By the local invariant cycle theorem [11],
$ \eta' $ comes from an element
$ \eta'_{0} $ of
$ H^{2}(X_{0},\bQ) $,
where
$ X_{s} = \{s z_{0}^{d} + z_{0}f_{d-1} + f_{d} = 0\}
\subset \bP^{3} $ (here
$ f_{0} $ is identified with
$ s) $.

If
$ \eta $ is viewed as a family of cohomology classes by
Poincar\'e duality, the integral
$ \int_{\eta}\omega_{s} $ is defined by using the pairing of
cohomology classes.
By (3.7) below,
$ \{\omega_{s}\} $ is extended to a section
$ \tom $ of
$ \cL $,
because letting
$ \tf = f/z_{0}^{d} $ and
$ x_{i} = z_{i}/z_{0} $, the restriction of
$ \omega_{s} $ to
$ X'_{s} := X_{s} \setminus z_{0}^{-1}(0) $ for
$ s \ne 0 $ is given by
$$
-\Res_{X'_{s}}\tf^{-1}x_{3}^{d-4}dx_{1}\wedge dx_{2}\wedge dx_{3}
= -(x_{3}^{d-4}dx_{1}\wedge dx_{2}\wedge dx_{3}/d\tf)|_{X'_{s}}.
$$
(Note that
$ m, n $ and
$ d $ in (3.7) are respectively
$ d - 4, 3 $ and
$ d - 1 $.)
Then the limit of the pairing of
$ \omega_{s} $ and
$ \eta $ for
$ s \to 0 $ depends only on the image of
$ \tom $ in
$ \Gr_{V}^{0}\cL(0) $ and the image of
$ \eta $ in
$ \Gr_{F}^{0}\Gr_{V}^{0}\cL(0) $.
This can be verified by expressing
$ \tom $ as a sum of
$ v(t)\tu $ where
$ v(t) $ is a holomorphic function of
$ t $, and
$ \tu $ is as in (3.4.2) below for
$ u \in L_{\infty,e(-a)} $ with
$ a \in [0,1) $, because the pairing is defined for
local systems.
So we may replace
$ \eta $ with
$ \eta' $ as long as we consider the limit of the integral for
$ s \to 0 $.

Let
$ B $ be a sufficiently small ball with center
$ O := (0, 0, 0,1) $ defined by using the coordinates
$ x_{i} $ as above.
Let
$ Y_{s} = X_{s} \setminus B $ with the inclusion
$ j_{s} : Y_{s} \to X_{s} $.
Since
$ X_{0} \cap B $ is contractible, we have a canonical
isomorphism
$ H_{2}(Y_{0},\bQ) = H_{c}^{2}(Y_{0},\bQ) =
H^{2}(X_{0},\bQ) $,
and
$ \eta'_{0} $ is identified with an element of
$ H_{2}(Y_{0},\bQ) = H_{c}^{2}(Y_{0},\bQ) $ where the
last isomorphism comes from Poincar\'e duality.
(Here we omit Tate twists to simplify the notation.)
Since
$ \{Y_{s}\}_{s\in \Delta} $ is a topologically trivial family,
$ \eta'_{0} $ is extended to a section of the constant local system
$ \{H_{2}(Y_{s},\bQ)\}_{s\in \Delta} $ or
$ \{H_{c}^{2}(Y_{s},\bQ)\}_{s\in \Delta} $ which is identified with
$ \eta' $ by the injection
$ H_{2}(Y_{s},\bQ) \to H_{2}(X_{s},\bQ) $ or
$ H_{c}^{2}(Y_{s},\bQ) \to H^{2}(X_{s},\bQ) $.
(Note that
$ H^{1}(X_{s} \cap B,\bQ) = 0 $ by the theory of Milnor fibration [30].)
Thus the integral
$ \int_{\eta'}\omega_{s} $ is defined by restricting
$ \omega_{s} $ to
$ Y_{s} $.
This is well-defined also for
$ s = 0 $.

Let
$ \tX_{0} $ be the blow-up of
$ X_{0} $ at
$ O $, and
$ C $ be the exceptional divisor.
Then
$ \tX_{0} $ is the blow-up of
$ \bP^{2} $ along
$ f_{d-1}^{-1}(0) \cap f_{d}^{-1}(0) $,
and
$ C $ is the proper transform of
$ f_{d-1}^{-1}(0) $,
because
$ f_{d-1}^{-1}(0) $ intersects
$ f_{d}^{-1}(0) $ transversally at smooth points.
Consider the exact sequence
$$
0 \to H^{1}(C,\bQ) \overset{\iota}\to\to H_{c}^{2}
(\tX_{0}\setminus C,\bQ) \to H^{2}(\tX_{0},\bQ),
$$
where
$ \iota $ is the dual of
$ \Res_{C} : H^{2}(\tX_{0}\setminus C,\bQ) \to H^{1}(C,\bQ)(-1) $.
Since the last term of the exact sequence is a direct sum of
$ \bQ(-1) $,
$ \iota $ induces an isomorphism
$$
H^{1}(C,\bC)/F^{1} = H_{c}^{2}(\tX_{0}\setminus C,\bC)/F^{1}.
\leqno(3.3.4)
$$
Let
$ \eta'' \in H^{1}(C,\bC)/F^{1} $ corresponding to
$ \eta'_{0} $ (mod
$ F^{1}) $.
This is independent of the choice of
$ \eta' $ by (3.5).
We have
$$
\langle \eta'',\Res_{C}\omega_{0} \rangle = \langle \eta'_{0},
\omega_{0}\rangle \,(= \int_{\eta'_{0}}\omega_{0}),
$$
where
$ \langle *,* \rangle $ denotes the scalar extension of the pairings
$$
\aligned
H^{1}(C,\bQ)\otimes H^{1}(C,\bQ)(-1)
&\to \bQ(-2),
\\
H_{c}^{2}(\tX_{0}\setminus C,\bQ)\otimes H^{2}(\tX_{0}\setminus C,\bQ)
&\to \bQ(-2).
\endaligned
$$
Note that the last pairing is a perfect pairing of mixed Hodge
structures and
$ \langle F^{1},\omega_{0} \rangle = 0 $ because
$ \omega_{0} \in F^{2}H^{2}(\tX_{0}\setminus C,\bC) $) (i.e.
$ \omega_{0} $ is a logarithmic
$ 2 $-form on
$ \tX_{0}\setminus C $).
Indeed, let
$ \omega' $ be the logarithmic
$ 2 $-from on
$ \bP^{2} \setminus f_{d-1}^{-1}(0) $ which is expressed as
$ dy_{1}\wedge dy_{2}/\of_{d-1} $ using
the coordinates
$ y_{i} $ on
$ \{z_{3} \ne 0\} $ as in (3.3.1).
Then
$ \omega_{0} $ is the pull-back of
$ \omega' $ by the blow-up of
$ \bP^{2} $ along
$ f_{d-1}^{-1}(0) \cap f_{d}^{-1}(0) $.
(Note that it does not have a pole along the exceptional divisor
of the blow-up, because the pole is cancelled by the zero coming
from the pull-back of a
$ 2 $-from.)
This shows also that
$ \Res_{C}\omega_{0} $ is independent of
$ g $ (using the isomorphism
$ C = f_{d-1}^{-1}(0)) $.
We will show that
$ \eta'' $ is also independent of
$ g $.

Since we can construct
$ \eta'_{0} $ depending continuously on
$ g $ locally on
$ S' $ by (3.6) below, it gives a local section of the local system
$ \{H_{c}^{2}(\tX_{0}\setminus C,\bQ)\} $ on
$ S' $, but this is not unique.
However, it induces a global section
$ \tet $ of a quotient local system
$ L' $ divided by a certain subsheaf which underlies a
variation of Hodge structure of type
$ (1.1) $.
(Indeed,
$ \eta'_{0} $ is unique up to
$ F^{1}H_{c}^{2}(\tX_{0}\setminus C,\bQ) $,
and if
$ g $ does not belong to the subspace of
$ S' $ on which
$ \dim F^{1}H_{c}^{2}(\tX_{0}\setminus C,\bQ) $ jumps, then
$ F^{1}H_{c}^{2}(\tX_{0}\setminus C,\bQ) $ underlies the stalk at
$ g $ of a variation of Hodge structure of type
$ (1,1) $ on
$ S' $.
Note that
$ \dim F^{1}H_{c}^{2}(\tX_{0}\setminus C,\bQ) $ is constant outside a
subset which is locally a countable union of proper analytic subsets.)

Let
$ L'' $ be the quotient local system of
$ L' $ divided by the constant local system
$ \{H^{1}(C,\bQ)\} $ (using the canonical isomorphism
$ C = f_{d-1}^{-1}(0) $).
Then
$ L', L'' $ underlie variations of Hodge structures
$ H', H'' $ on
$ S' $, and
$ \tet $ induces a global section of
$ L'' $ which gives a morphism of Hodge structures
$ \bQ(-1) \to H'' $, because
$ H'' $ is of type
$ (1,1) $.
Taking the pull-back of the short exact sequence
$$
0 \to \{H^{1}(C,\bQ)\} \to H' \to H'' \to 0
$$
by this morphism, we get an extension of
$ \bQ(-1) $ by the constant variation
$ \{H^{1}(C,\bZ)\} $ in the category of admissible variations of
mixed Hodge structures on
$ S' $.

Let
$ \tau_{g} $ denote the image of
$ \eta'' $ in
$$
J(C) = H^{1}(C,\bC)/(F^{1} + H^{1}(C,\bZ)(1)).
$$
Then, using [10] and an isomorphism similar to (3.3.4) (with
$ H_{c}^{2}(\tX_{0}\setminus C,\bC) $ replaced by a stalk of
$ H' $), we see that
$ \{\tau_{g}\} $ is the admissible normal function corresponding to
the above extension.
(This argument is inspired by Deligne's reformulation of Griffiths'
Abel-Jacobi map, see [20].)
Furthermore
$ \{\tau_{g}\} $ is constant by (2.1.1) and (2.1.2),
because
$ H := \{H^{1}(C,\bQ)(1)\} $ is a constant variation of Hodge
structure of weight
$ -1 $ on a smooth affine rational curve
$ S' $ so that
$ \Gr_{0}^{W}H^{1}(S',H) = 0 $.
Thus
$ \eta'' $ is independent of
$ g $, because it depends on
$ g $ continuously.

On the other hand,
$ \Res_{C}\omega_{0} $ is also independent of
$ g $ as seen above.
Therefore
$ \langle \eta'',\Res_{C}\omega_{0} \rangle $ and hence
$ \int_{\eta'_{0}}\omega_{0} $ are independent of
$ g $.
But this integral coincides with the limit of
$ \int_{\Gamma}\omega_{s} $ for
$ s \to 0 $ by the above argument.
The latter is equal to
$ \int_{\Gamma}\omega_{0} $,
which is a nonconstant function of
$ g $ as shown above.
This is a contradiction, and the assertion follows.

\medskip
To complete the proof of Theorem (3.3), we need some knowledge
of Deligne extension [15] and limit mixed Hodge
structure ([39], [42]):

\medskip\noindent
{\bf 3.4.~Complement to the proof of (3.3),
I.~Limit mixed Hodge structure.}
Let
$ L $ be a local system with rational coefficients on a
punctured disk
$ \Delta^{*} $ such that its monodromy
$ T $ is quasi-unipotent.
Let
$ \cL $ the Deligne extension of
$ \cL^{*} := \cO_{\Delta^{*}}\otimes_{\bQ}L $ such that
the eigenvalues of the residue of the connection at the origin
are contained in
$ [0,1) $,
see [15].
Let
$ p : \tDe^{*} \to \Delta^{*} $ be a universal
cover, and define
$ L_{\infty} = \Gamma (\tDe^{*},p^{*}L) $.
Let
$ \cL(0) $ denote the fiber of
$ \cL $ at the origin (i.e.
$ \cL\otimes_{\cO_{\Delta}}(\cO_{\Delta}/m_{\Delta,0}) $
where
$ m_{\Delta,0} $ is the maximal ideal at the origin).
Then, choosing a coordinate
$ t $ of
$ \Delta $,
we have a canonical isomorphism
$ L_{\infty} = \cL(0) $ induced by
$$
u \in L_{\infty} \mapsto \tu :=
\exp\biggl(-\frac{\log t}{2\pi i}\log T\biggr)u
\in \Gamma (\Delta,\cL),
\leqno(3.4.1)
$$
where
$ \log T $ is the logarithm of the monodromy
$ T $ whose eigenvalues divided by
$ -2\pi i $ are contained in
$ [0,1) $.
Let
$ T = T_{s}T_{u} $ be the Jordan decomposition, and put
$ N = \log T_{u} $.
Let
$ L_{\infty,\lambda} = \Ker(T_{s} - \lambda) \subset L_{\infty} $
for
$ \lambda \in \bC $, and
$ e(\alpha) = \exp(2\pi i\alpha) $ for
$ \alpha \in \bQ $.
Then for
$ u \in L_{\infty,e(-\alpha)} $ with
$ \alpha \in [0,1) $, we have
$$
\tu = \sum_{i\ge 0}t^{\alpha}\biggl(-\frac{\log t}{2\pi i}\biggr)^{i}
N^{i}u/i!
\leqno(3.4.2)
$$

Assume that
$ L $ underlies a polarizable variation of Hodge structure.
Then the Hodge filtration
$ F $ on
$ \cL^{*} $ can be extended to that of
$ \cL $,
see [39].
Let
$ F $ denote also the quotient filtration on
$ \cL(0) $.
It does not necessarily give the the Hodge filtration of
the limit mixed Hodge structure unless
$ T $ is unipotent.
We have to take further the graded pieces of the filtration
$ V $ on
$ \cL $ where
$ V^{\alpha} $ is the Deligne extension such that the eigenvalues
of the residue of the connection are contained in
$ [\alpha,\alpha+1) $ for
$ \alpha \ge 0 $.
(This filtration comes essentially from the
$ m $-adic filtration on a ramified base change of
$ \Delta $ such that the pull-back of
$ L $ has unipotent monodromy as in loc.~cit).
It induces a decreasing filtration
$ V $ on
$ L_{\infty} $ such that
$$
V^{\alpha}L_{\infty} = \oplus_{\alpha \le \beta <1}
L_{\infty,e(-\beta)}\quad \text{for}\,\,\alpha \in [0,1].
$$
Then the limit Hodge filtration is given by
$ \Gr_{V}F $.
(This is clear if we consider a unipotent base change as above.)

\medskip\noindent
{\bf 3.5.~Lemma.}
{\it With the above notation, assume the variation of Hodge
structure has weight
$ 2 $ and level
$ 2 $ {\rm (}i.e.
$ \Gr_{F}^{p}\cL = 0 $ unless
$ p \in [0,2]) $.
Let
$ \eta $ be a multivalued section of
$ L $ such that the variation
$ T\eta - \eta $ belongs to
$ F^{1}\cL^{*} $.
Then
$ \eta $ is extended to a section of
$ \Gr_{F}^{0}\cL $ and its image in
$ \Gr_{F}^{0}\Gr_{V}^{0}\cL(0) $ coincides with the image of
some
$ \eta' \in \Gamma(\Delta^{*},L) $ which is unique up to a
section of
$ F^{1}\cL^{*} $.
}

\medskip\noindent
{\it Proof.}
Let
$ L_{1} $ be the subsheaf of
$ L $ generated by
$\bQ[T](T\eta-\eta) $.
Then it underlies a variation of Hodge structure of type (1,1),
and is a direct factor of
$ L $ by the semisimplicity of a polarizable variation of
Hodge structure.
So we may assume either
$ L_{1} = 0 $ or
$ L = L_{1} $.
Then the assertion is clear.

\medskip\noindent
{\bf 3.6.~Remark.}
Let
$ S' $ be an open disk, and
$ L $ be a local system on
$ S'\times\Delta^{*} $ which has quasi-unipotent monodromy,
and underlies a variation of Hodge structure as in (3.5).
Then the Deligne extension
$ \cL $ on
$ S'\times\Delta $ and the filtration
$ V $ are defined as above, and the assertion of
(3.5) holds in this setting.

\medskip\noindent
{\bf 3.7.~Complement to the proof of (3.3),
II.~Relative logarithmic forms.}
Let
$ f : Y \to S $ be a projective morphism of complex manifolds of
relative dimension
$ n $,
where
$ S $ is an open disk with a coordinate
$ t $.
Assume
$ f $ is smooth over the punctured disk
$ S^{*} $,
and the singular fiber
$ D $ is a divisor with normal crossings (but not necessarily
reduced).
Let
$ {\Omega}_{Y}^{\ssbul}(\log D) $ be the complex of logarithmic
forms.
Then the Koszul complex
$ K^{\ssbul}({\Omega}_{Y}^{\ssbul}(\log D), f^{*}(dt/t)\wedge ) $
is acyclic, and the relative logarithmic forms
$ {\Omega}_{Y/S}^{\ssbul}(\log D) $ are defined to be the
cokernel of
$ f^{*}(dt/t)\wedge $,
see [42].
It is well-known that the higher direct images
$ R^{j}f_{*}{\Omega}_{Y/S}^{\ssbul}(\log D) $ are locally free
$ \cO_{S} $-Modules with a logarithmic connection
$ \nabla $ such that the eigenvalues of the residue of
$ \nabla_{t\partial /\partial t} $ are contained in
$ [0,1) $,
and it gives the Deligne extension of
$ (R^{j}f_{*}\bQ_{X}|_{S^{*}})\otimes_{\bQ}\cO_{S^{*}} $,
see loc.~cit.

Let now
$ f' : X \to S $ be a projective morphism of complex manifolds
such that
$ \Sing f' = \{O\} $,
where
$ S $ is as above.
Assume that the inverse image of the singular fiber
$ X_{0} $ by the blowing-up
$ \pi : Y \to X $ with center
$ O $ is a divisor with normal crossings.
Then the intersection of the proper transform of
$ X_{0} $ with the exceptional divisor is a smooth hypersurface,
and its degree
$ d $ coincides with the multiplicity of
$ X_{0} $ at
$ O $.
Let
$ \omega $ be a holomorphic
$ (n+1) $-form on
$ X $,
and
$ m $ be the multiplicity of the zero of
$ \omega $ at
$ O $.
Then
$ \pi^{*}\omega $ and
$ (f'\pi )^{*}t $ has zeros of order
$ m + n - 1 $ and
$ d $ respectively at the exceptional divisor.
So, if
$ m + n - 1 \ge d - 1 $,
then
$ \pi^{*}\omega /(f'\pi )^{*}t $ is a logarithmic
$ (n+1) $-form on
$ X $,
and is identified with a relative logarithmic
$ n $-form by the acyclicity of the above Koszul complex.

\medskip\noindent
{\bf 3.8.~Remark.}
It is not easy to construct an indecomposable higher cycle on a
given variety.
Any known analytic method constructs a family of varieties with a
higher cycle, and shows for a generic member that the cycle is
indecomposable, or it is not annihilated by the reduced
Abel-Jacobi map, see [12], [13], [14], [31], [36], [45].
The situation is different over a number field (see [29]),
although it is not easy to express a cycle explicitly.
Note that once we get a variety and a higher cycle satisfying the
last condition in Theorem (0.1), there is a finitely generated
subfield
$ k $ over which the variety and the cycle are defined.

\bigskip\bigskip
\centerline{{\bf 4.~Inductive construction}}

\bigskip\noindent
{\bf 4.1.~Borel-Moore cohomology.}
For a singular variety
$ Z $,
let
$ \bD_{Z} = a_{Z}^{!}\bQ \in D_{c}^{b}(Z,\bQ) $ with
$ a_{Z} : Z \to pt $ the structure morphism.
Here
$ D_{c}^{b}(Z,\bQ) $ is the full subcategory of
$ D_{c}^{b}(Z^{\an},\bQ) $ consisting of complexes with
algebraic stratifications.
Then the Borel-Moore homology
$ H_{j}^{\BM}(Z,\bQ) $ is given by
$ H^{-j}(a_{Z})_{*}\bD_{Z} $.
We define the Borel-Moore cohomology by
$ H_{\BM}^{j}(Z,\bQ) = H_{2n-j}^{\BM}(Z,\bQ)(-n) $ if
$ Z $ is purely
$ n $-dimensional.
If
$ Z $ is smooth, we have
$ \bD_{Z} = \bQ_{Z}(n)[2n] $,
and
$ H_{\BM}^{j}(Z,\bQ) = H^{j}(Z,\bQ) $.

\medskip\noindent
{\bf 4.2.~Construction.}
Let
$ f \in \bC[z_{0}, \dots, z_{n}], h \in \bC[w_{0}, w_{1}] $
be homogeneous polynomials of degree
$ d $,
where
$ n \ge 1 $,
$ d \ge 2 $.
Put
$ \tf = f + h, X = f^{-1}(0) \subset \bP^{n} $,
$ Z = h^{-1}(0) \subset \bP^{1} $,
$ \tX = \tf^{-1}(0) \subset \bP^{n+2} $.
Assume
$ X, Z $ are smooth (in particular,
$ Z $ consists of
$ d $ distinct points in
$ \bP^{1}) $.
Then
$ \tX $ is also smooth.
Let
$ p_{i} = (a_{i}, b_{i}) $ be the points of
$ Z $ for
$ 1 \le i \le d $.
Then the equation
$ b_{i}w_{0} = a_{i}w_{1} $ defines a subvariety
$ Y_{i} $ of
$ \tX $.
It is a projective cone over
$ X $,
and has a unique singular point
$ q_{i} $.
Let
$ \tY_{i} $ be the blow-up of
$ Y_{i} $ at
$ q_{i} $.
It is a
$ \bP^{1} $-bundle over
$ X $.
Let
$ \pi_{i} : \tY_{i} \to X $ and
$ \rho_{i} : \tY_{i} \to \tX $ denote the
canonical morphisms.
Then we get a morphism of higher Chow groups
$$
(\rho_{i})_{*}\pi_{i}^{*} : \CH^{p}(X,m) \to \CH^{p+1}(\tX,m)
\leqno(4.2.1)
$$
which is compatible with the morphism in cohomology
$$
(\rho_{i})_{*}\pi_{i}^{*} : H^{j}(X,\bQ) \to H^{j+2}(\tX,\bQ)(1)
\leqno(4.2.2)
$$
via the cycle map (1.1.1).
Here we assume
$ m > 0 $ so that the target of the cycle map is
$$
\Ext_{\MHS}^{1}(\bQ,H^{2p-m-1}(X,\bQ)(p))
\leqno(4.2.3)
$$
and similarly for
$ \tX $.
Since
$ H^{j}(X,\bQ) $ is interesting only for
$ j = n - 1 $,
we may assume
$ 2p = n + m $.

\medskip\noindent
{\bf 4.3.~Proposition.}
{\it The morphism {\rm (4.2.2)} is injective for
$ j = n - 1 $.
}

\medskip\noindent
{\it Proof.}
Let
$ Y'_{i} = Y_{i} \setminus \{q_{i}\} $.
It is a line bundle over
$ X $,
and the zero section is the intersection with another
$ Y_{j} $.
So we get for
$ i \ne j $ a cartesian diagram
$$
\CD
\tX @<{\rho_{j}}<< Y_{j}
\\
@AA{\rho_{i}}A @AA{\iota_{j}}A
\\
Y_{i} @<{\iota_{i}}<< X
\endCD
$$
where
$ \iota_{i} : X \to Y_{i} $ denotes the zero section.
We have
$ \bD_{Y'_{i}} = \bQ_{Y'_{i}}(n)[2n] $,
and hence
$ \iota_{i}^{*}\bD_{Y_{i}} = \bQ_{X}(n)[2n] $.
So with the notation of (4.1), we get the restriction morphism
$$
\iota_{i}^{*} : H_{\BM}^{n-1}(Y_{i},\bQ) \to H^{n-1}(X,\bQ).
\leqno(4.3.1)
$$
This is an isomorphism because the restriction morphism
$$
H_{\BM}^{n-1}(Y_{i},\bQ) \to H^{n-1}(Y'_{i},\bQ)
$$
is an isomorphism by the long exact sequence of local cohomology
combined with the isomorphism
$ \iota_{q_{i}}^{!}\bD_{Y_{i}} = \bQ $, where
$ \iota_{q_{i}} : \{ q_{i} \} \to Y_{i} $ denotes the inclusion.
By duality, we get the isomorphism
$$
(\iota_{j})_{*} : H^{n-1}(X,\bQ) \to H^{n+1}(Y_{j},\bQ)(1).
\leqno(4.3.2)
$$

Since
$ \rho_{j}^{*}\bQ_{\tX} = \bQ_{Y_{i}} $,
the canonical morphism
$ (\rho_{i})_{*} : \bD_{Y_{i}} \to \bD_{\tX} $
together with the functorial morphism
$ id \to (\rho_{j})_{*}\rho_{j}^{*} $ induces a commutative diagram
$$
\CD
\bD_{\tX} @>{\rho_{j}^{*}}>> \bQ_{Y_{j}}(n+1)[2n+2]
\\
@AA{(\rho_{i})_{*}}A @AA{(\iota_{j})_{*}}A
\\
\bD_{Y_{i}} @>{\iota_{i}^{*}}>> \bQ_{X}(n)[2n]
\endCD
$$
where the direct images by closed embeddings are omitted to
simplify the notation.
So we get
$ \rho_{j}^{*}(\rho_{i})_{*} = (\iota_{j})_{*}\iota_{i}^{*} :
H_{\BM}^{n-1}(Y_{i},\bQ) \to H^{n+1}(Y_{j},\bQ)(1) $.
Since the second composition is an isomorphism by the above argument,
we get the injectivity of
$ (\rho_{i})_{*} : H_{\BM}^{n-1}(Y_{i},\bQ) \to H^{n+1}
(\tX,\bQ)(1) $.

Now it remains to show the injectivity of the composition of
the restriction morphism
$ \pi_{i}^{*} : H^{n-1}(X,\bQ) \to H^{n-1}(\tY_{i},\bQ) $ and the
Gysin morphism
$ H^{n-1}(\tY_{i},\bQ) \to H_{\BM}^{n-1}(Y_{i},\bQ) $.
But its further composition with the restriction morphism to
$ H^{n-1}(Y'_{i},\bQ) $ is the restriction morphism
$ (\pi'_{i})^{*} : H^{n-1}(X,\bQ) \to H^{n-1}(Y'_{i},\bQ) $
where
$ \pi'_{i} : Y'_{i} \to X $ denotes the projection.
This is clearly an isomorphism, and the assertion follows.

\medskip\noindent
{\bf 4.4.~Theorem.}
{\it With the notation of {\rm (4.2),} let
$ \zeta \in \CH^{p}(X,m) $ where
$ 2p = m + n $.
We consider the cycle map to {\rm (1.1.3)} if
$ m > 1 $ and to {\rm (1.1.4)} if
$ m = 1 $.
Assume
$ \zeta $ is not annihilated by this cycle map.
Then the image of
$ (\rho_{i})_{*}\pi_{i}^{*}\zeta \in \CH^{p+1}(\tX,m) $ by the
same cycle map for
$ \tX $ does not vanish.
In particular,
$ (\rho_{i})_{*}\pi_{i}^{*}\zeta $ is indecomposable if
$ m = 1 $.
}

\medskip\noindent
{\it Proof.}
This follows from (4.3) combined with the semisimplicity of
$ H^{n+1}(\tX,\bQ) $ as a Hodge structure.

\newpage\centerline{{\bf References}}

\bigskip
\item{[1]}
Asakura, M., Motives and algebraic de Rham cohomology, in: The
arithmetic and geometry of algebraic cycles (Banff), CRM Proc.
Lect. Notes, 24, AMS, 2000, pp. 133--154.

\item{[2]}
Beilinson, A., Higher regulators and values of
$ L $-functions, J. Soviet Math. 30 (1985), 2036--2070.

\item{[3]}
\SameAuthor, Notes on absolute Hodge cohomology, Contemporary Math. 55
(1986) 35--68.

\item{[4]}
\SameAuthor, Height pairing between algebraic cycles, Lect. Notes in
Math., vol. 1289, Springer, Berlin, 1987, pp. 1--26.

\item{[5]}
Beilinson, A., Bernstein, J. and Deligne, P., Faisceaux pervers,
Ast\'erisque, vol. 100, Soc. Math. France, Paris, 1982.

\item{[6]}
Bloch, S., Lectures on algebraic cycles, Duke University
Mathematical series 4, Durham, 1980.

\item{[7]}
\SameAuthor, Algebraic cycles and values of
$ L $-functions, J. Reine Angew. Math. 350 (1984), 94--108.

\item{[8]}
\SameAuthor, Algebraic cycles and higher
$ K $-theory, Advances in Math., 61 (1986), 267--304.

\item{[9]}
\SameAuthor, Algebraic cycles and the Beilinson conjectures,
Contemporary Math. 58 (1) (1986), 65--79.

\item{[10]}
Carlson, J., Extensions of mixed Hodge structures, in:
Journ\'ees de G\'eom\'etrie Alg\'ebrique d'Angers 1979,
Sijthoff-Noordhoff Alphen a/d Rijn, 1980, pp. 107--128.

\item{[11]}
Clemens, C.H., Degeneration of K\"ahler manifolds, Duke Math. J.
44 (1977), 215--290.

\item{[12]}
Collino, A., Griffiths' infinitesimal invariant and higher
$ K $-theory on hyperelliptic Jacobians, J. Alg. Geom. 6 (1997),
393--415.

\item{[13]}
\SameAuthor, Indecomposable motivic cohomology classes on quartic
surfaces and on cubic fourfolds, in Algebraic
$ K $-theory and Application (Eds H. Bass et al.),
World Scientific, 1999, pp. 370--402.

\item{[14]}
del Angel, P. and M\"uller-Stach, S.,
The transcendental part of the regulator map for
$ K_{1} $ on a mirror family of
$ K3 $ surfaces, preprint.

\item{[15]}
Deligne, P., Equations diff\'erentielles \`a points singuliers
r\'eguliers, Lect. Notes in Math., Vol. 163. Springer, Berlin, 1970.

\item{[16]}
\SameAuthor, Th\'eorie de Hodge I, Actes Congr\`es Intern. Math.,
1970, vol. 1, 425-430; II, Publ. Math. IHES, 40 (1971), 5--57;
III ibid., 44 (1974), 5--77.

\item{[17]}
\SameAuthor, Valeurs de fonctions
$ L $ et p\'eriodes d'int\'egrales,
Proc. Symp. in pure Math., 33 (1979) part 2, pp. 313--346.

\item{[18]}
Deligne, P., Milne, J., Ogus, A. and Shih, K., Hodge Cycles, Motives, and
Shimura varieties, Lect. Notes in Math., vol 900, Springer, Berlin, 1982.

\item{[19]}
Deninger, C. and Scholl, A., The Beilinson conjectures,
Proceedings Cambridge Math. Soc. (eds. Coats and Taylor) 153
(1992), 173--209.

\item{[20]}
El Zein, F. and Zucker, S., Extendability of normal functions
associated to algebraic cycles, in: Topics in transcendental
algebraic geometry, Ann. Math. Stud., 106, Princeton Univ. Press,
Princeton, N.J., 1984, pp. 269--288.

\item{[21]}
Esnault, H. and Viehweg, E., Deligne-Beilinson cohomology, in:
Beilinson's conjectures on Special Values of L-functions,
Academic Press, Boston, 1988, pp. 43--92.

\item{[22]}
Green, M., Griffiths' infinitesimal invariant and the Abel-Jacobi
map, J. Diff. Geom. 29 (1989), 545--555.

\item{[23]}
Griffiths, P., On the period of certain rational integrals I,
II, Ann. Math. 90 (1969), 460--541.

\item{[24]}
Jannsen, U., Deligne homology, Hodge-$ D $-conjecture, and motives,
in Beilinson's conjectures on Special Values of
$ L $-functions, Academic Press, Boston, 1988, pp. 305--372.

\item{[25]}
\SameAuthor, Mixed motives and algebraic
$ K $-theory, Lect. Notes in Math., vol. 1400, Springer, Berlin, 1990.

\item{[26]}
\SameAuthor, Motivic sheaves and filtrations on Chow groups, Proc. Symp.
Pure Math. 55 (1994), Part 1, pp. 245--302.

\item{[27]}
Lang, S., Abelian varieties, Interscience Publishers, New York,
1959.

\item{[28]}
Levine, M., Localization on singular varieties, Inv. Math. 91 (1988),
423--464.

\item{[29]}
Mildenhall, S.J.M., Cycles in a product of elliptic curves, and a group
analogous to the class group, Duke Math. J. 67 (1992), 387--406.

\item{[30]}
Milnor, J., Singular points of complex hypersurfaces, Ann. Math. Stud.
vol. 61, Princeton Univ. Press, 1969.

\item{[31]}
M\"uller-Stach, S., Constructing indecomposable motivic cohomology
classes on algebraic surfaces, J. Alg. Geom. 6 (1997), 513--543.

\item{[32]}
Rosenschon, A. and Saito, M., Cycle map for strictly
decomposable cycles, preprint (math.AG/0104079).

\item{[33]}
Saito, M., Mixed Hodge Modules, Publ. RIMS, Kyoto Univ., 26
(1990), 221--333.

\item{[34]}
\SameAuthor, Extension of mixed Hodge Modules, Compos. Math. 74 (1990),
209--234.

\item{[35]}
\SameAuthor, Admissible normal functions, J. Alg. Geom. 5 (1996),
235--276.

\item{[36]}
\SameAuthor, Bloch's conjecture, Deligne cohomology and higher
Chow groups, preprint RIMS--1284 (or math.AG/9910113).

\item{[37]}
\SameAuthor, Arithmetic mixed sheaves, Inv. Math. 144 (2001),
533--569.

\item{[38]}
\SameAuthor, Refined cycle maps, preprint (math.AG/0103116).

\item{[39]}
Schmid, W., Variation of Hodge structure: the singularities of
the period mapping, Inv. Math. 22 (1973), 211--319.

\item{[40]}
Schoen, C., Zero cycles modulo rational equivalence for some
varieties over fields of transcendence degree one, Proc. Symp.
Pure Math. 46 (1987), part 2, pp. 463--473.

\item{[41]}
\SameAuthor, On certain exterior product maps of Chow groups,
Math. Res. Let. 7 (2000), 177--194,

\item{[42]}
Steenbrink, J.H.M., Limits of Hodge structures, Inv. Math. 31
(1975/76), 229--257.

\item{[43]}
Voisin, C., Variations de structures de Hodge et z\'ero-cycles sur
les surfaces g\'en\'erales, Math. Ann. 299 (1994), 77--103.

\item{[44]}
\SameAuthor, Transcendental methods in the study of algebraic
cycles, Lect. Notes in Math. vol. 1594, pp. 153--222.

\item{[45]}
\SameAuthor, Nori's connectivity theorem and higher Chow groups, preprint.

\item{[46]}
Zucker, S., Hodge theory with degenerating coefficients,
$ L_{2} $-cohomology in the Poincar\'e metric,
Ann. Math., 109 (1979), 415--476.

\bigskip\noindent
Andreas Rosenschon

\noindent
Department of Mathematics, Duke University,
Durham, NC 27708, U.S.A

\noindent
E-Mail: axr\@math.duke.edu

\medskip\noindent
Morihiko Saito

\noindent
RIMS Kyoto University,
Kyoto 606--8502 Japan

\noindent
E-Mail: msaito\@kurims.kyoto-u.ac.jp

\bigskip\noindent
\ver

\bye